\documentclass[11pt]{article}
\usepackage{amssymb}
\usepackage{amscd}
\usepackage{epsf}
\newtheorem{thm}{Theorem}[section]

\newtheorem{prop}[thm]{Proposition}

\newtheorem{defi}[thm]{Definition}

\newtheorem{ex}[thm]{Example}
\newcommand{\inv}{^{-1}}

\newcommand{\C}{{\mathbb{C}}}
\newcommand{\Z}{{\mathbb{Z}}}
\newcommand{\R}{{\mathbb{R}}}

\newcommand{\T}{{\mathbb{T}}}

\newcommand{\Gm}{{\mathbb{G}}_m}

\newcommand{\ov}{\overline}

\newcommand{\Hom}{\mbox{\rm Hom}}
\newcommand{\diag}{\mbox{\rm diag}}
\newcommand{\trop}{\mbox{\rm trop}}

\newcommand{\Rinf}{\R \cup \{-\infty\}}
\begin{document}

\title{A tropical view on Bruhat-Tits buildings and their compactifications}
\author{Annette Werner}
\date{}

\maketitle

\centerline{\bf Abstract:} We relate some features of Bruhat-Tits buildings and their compactifications to tropical geometry. 
If $G$ is a semisimple group over a suitable non-Archimedean field, the stabilizers of points in the Bruhat-Tits building of $G$ and in some of its compactifications are described by tropical linear algebra. The compactifications we consider arise from algebraic representations of $G$. We show that the fan which is used to compactify an apartment in this theory is given by the weight polytope of the representation and that it is related to the tropicalization of the hypersurface given by the character of the representation.
\small 
\\[0.3cm]

\centerline{{\bf 2010 MSC: }20E42, 20G25, 14T05 } 
\centerline{{\bf Keywords: }Bruhat-Tits buildings, tropical geometry, weight polytopes}

\normalsize 
\section*{Introduction}

Let $K$ be a field which is complete with respect to  a non-trivial discrete valuation and has perfect  residue field, and let $G$ be a semisimple group over $K$.  Then the Bruhat-Tits building $\mathfrak{B}(G)$ is a metric space with a continuous $G(K)$-action and a polysimplicial structure. It consists of apartments, which can be identified with real cocharacter spaces of maximal split tori. 
In the case of the group $SL_{n,K}$ we endow such an apartment with tropical addition and tropical multiplication. We say that a matrix $g$ in $SL_n(K)$ stabilizes a point $x$ in the apartment tropically if the real matrix we get by applying minus the valuation to $g$
stabilizes $x$ under tropical matrix multiplication. 

For every faithful representation $\rho: G \rightarrow SL_{n,K}$ there exists a $G(K)$-equivariant embedding $i: \mathfrak{B}(G) \rightarrow \mathfrak{B}(SL_{n,K})$ of Bruhat-Tits buildings, where $G(K)$ acts via $\rho$ on the right hand side. Our main result says that for every $x$ in $\mathfrak{B}(G)$ the subgroup $P_x$ of elements in $G(K)$ stabilizing $x$ coincides with the set of elements $g \in G(K)$ such that $\rho(g)$ stabilizes $i(x)$ tropically. Hence although  the group $P_x$ involves the action of $G(K)$ on the whole building its tropical interpretation takes place in one apartment. 

In order to present examples of subgroups of classical groups which are tropical stabilizers by our theorem, we recall the fact that in a simply connected group $G$ the stabilizer groups $P_x$ are the parahoric subgroups of $G(K)$. Hence they can be described as groups of matrices whose reduction modulo the valuation ideal lies in a certain parabolic subgroup over the residue field. We give a detailled treatment of the examples $G = SL_{n,K}$ and $G = Sp_{2n,K}$.  

Moreover we discuss compactifications of buildings associated to representations $\rho: G \rightarrow SL_{n,K}$ as in \cite{we3}. Every apartment is compactified with the help of a fan $\mathcal{F}_\rho$ defined by the combinatorics of the weights of $\rho$. We show that $\mathcal{F}_\rho$  is in fact the normal fan of the weight polytope of $\rho$. Besides, we show that the codimension one skeleton of $\mathcal{F}_\rho$ is the tropicalization of a hypersurface given by the character of $\rho$. If $G = SL_{n,K}$  and $K$ has characteristic zero, then this character is a Schur polynomial. We also prove a generalization of our main theorem for points in the boundary of the compactified building. 

Note that in the paper \cite{jsy} another tropical view on convexity in the building for $SL_n$ is presented.

The present paper is organized as follows. Section one provides the tropical linear algebra setting. Section two starts with the necessary facts about Bruhat-Tits buildings. In Proposition 2.4 the tropical interpretation of the stabilizer groups is given in the case $G = SL_{n,K}$.  Theorem 2.5 contains the general case. Section three deals with compactifications of Bruhat-Tits buildings. Propositions 3.5 and 3.6 relate the fan $\mathcal{F}_\rho$ to the weight polytope and to the tropicalization of the character of $\rho$. Theorem 3.9 contains the tropical interpretation of stabilizer groups for boundary points.

{\bf Acknowledgements: }I thank Michael Joswig, Bernd Sturmfels and Thorsten Theobald for useful discussions on tropical geometry. I am also grateful to MSRI for its hospitality during the program on Tropical Geometry, where  parts of this paper were written.

\section{The tropical torus} Tropical geometry is based on the tropical semiring $(\mathbb{R}, \oplus, \odot)$ with  $a \oplus b = \max\{a,b\}$ and $a \odot b = a +b$. Some authors use  the $(\min, +)$-version of the tropical semiring instead of the $(\max,+)$ version.

The space $\R^n$ together with componentwise addition $\oplus$ is a semimodule under the semiring $(\R, \oplus, \odot)$, if we put $a \odot (x_1,\ldots, x_n) = (a+ x_1, \ldots, a+ x_n)$ for $a \in \R$ and $(x_1,\ldots, x_n) \in \R^n$.  Let
\[ \T^{n-1} = \R^n / \R (1,\ldots, 1)\]
be the quotient
of $\R^n$ after the following equivalence relation:  $(x_1,\ldots, x_n) \sim (y_1,\ldots, y_n)$ if and only if there exists some $a \in \R$ such that $x_i = a \odot y_i = a + y_i $ for all $i$. We call $\T^{n-1}$ the tropical torus of rank $n-1$ as in \cite{jo}. We endow $\T^{n-1}$ with the quotient topology.

Let $K$ be a field with a non-trivial valuation map $v: K^\ast \rightarrow \R$. We put $v(0) = \infty$ and consider
the  negative valuation  map $-v: K \rightarrow \mathbb{R}_{- \infty} = \mathbb{R} \cup \{-\infty\}$. 
We extend the operations $\odot$ and $\oplus$ to $\R_{-\infty}$ by  $a \odot - \infty = - \infty$ for all $a \in \R_{-\infty}$ and $a \oplus -\infty = a$ for all $a \in \R_{-\infty}$. 

\begin{defi}
Let $g = (g_{ij})_{i,j}$ be a matrix in $GL_n(K)$. Then we define the associated tropical matrix as
$g_{\trop} = (-v ( g_{ij}))_{i,j} \in \mbox{Mat}_{n \times n} (\R)$. For every vector $x = ~^t(x_1,\ldots, x_n) \in \R^n$ we define the 
  vector $ g_{\trop} \cdot x = ~^t(y_1, \ldots, y_n) \in \R^n$ by the tropicalized linear action
 \[ y_i = -v(g_{i1}) \odot x_1 \oplus \ldots \oplus -v(g_{in}) \odot x_n = \max_j \{ -v(g_{ij}) + x_j\}.\]
\end{defi}

Note that $y_i$ lies indeed in $\R$, since at least one entry $g_{ij}$ in the $i$-th line must be non-zero, so that at least for one $j$ the term $-v(g_{ij}) + x_j$ is not equal to $-\infty$. 

Note that this does not define an action of $GL_n(K)$ on $\R^n$, as the following example shows. Take $n = 2$, and put
\[g = \left( \begin{array}{ll} 1 & 1\\ 0 & 1 \end{array} \right) \mbox{ and } h = \left( \begin{array}{ll} 1 & 0 \\ -1 & 1 \end{array} \right). \mbox{ Then}\]
\[(gh)_{\trop} \cdot \left(\begin{array}{l} x_1 \\ x_2 \end{array} \right) = \left(\begin{array}{l} x_2 \\ \max\{x_1, x_2\} \end{array} \right), \mbox{ but}\]
\[g_{\trop} \cdot \left( h_{\trop} \cdot \left(\begin{array}{l} x_1 \\ x_2 \end{array} \right)\right) = \left(\begin{array}{l} \max\{x_1, x_2\} \\ \max\{x_1, x_2\}\end{array} \right).\] 

The tropical matrix action describes a linear map in the context of Max-Plus Algebra, see \cite{abg} for an overview of this topic and further references. 

Let $T$ be a split torus over $K$ with character group $X^\ast(T) = \mbox{Hom}_K (T, \mathbb{G}_m)$ and cocharacter group $X_\ast(T) = \mbox{Hom}_K(\mathbb{G}_m, T)$. Let $f \in K[X^\ast(T)] = \Gamma(T, \mathcal{O}_T)$ be a non-zero global section of $T$. Hence $f$ is a Laurent polynomial in any basis $z_1,\ldots, z_d$ of $X^\ast(T)$. By the natural perfect pairing 
\[X^\ast(T) \times X_\ast(T) \rightarrow \Z\]
between characters and cocharacters, we evaluate characters on elements in $X_\ast(T)_\R = X_\ast(T) \otimes_\Z \R$. If $f = \sum_M a_M M$ with  $M \in X^\ast(T)$ and $a_M \in K$, then the tropical hypersurface associated to $f$ is defined as 
\[\mathcal{T}(f) = \{ x \in X_\ast(T)_\R: 
 \max_{M} \{ -v(a_M) + M(x) \} \mbox{ is attained at least twice}\}. 
\]
Choose a basis $z_1, \ldots, z_d$ of $X^\ast(T)$, i.e. an isomorphism $\varphi: \mathbb{G}_m^d \rightarrow T$. Then $f$ can be written as $f = \sum_I a_I z_1^{i_1} \ldots z_d^{i_d}$. The corresponding isomorphism $X_\ast(T)_\R \simeq X_\ast(\mathbb{G}_m^d)_\R = \R^d$ maps  $\mathcal{T}(f)$ bijectively to the subset 
\[\{ (u_1, \ldots, u_d) \in \R^d:   \max_I \{-v(a_I) + i_1 u_1 +  \ldots + i_d u_d\}\mbox{ is attained at least twice} \}. \]
By \cite{ekl}, Theorem 2.1.1, this set is equal to  the closure of the subset
\[\{ (-v(\alpha_1), \ldots, -v(\alpha_d)) \in \R^d: \alpha_1, \ldots, \alpha_d \in \overline{K}^\times \mbox{ with }f(\alpha_1, \ldots, \alpha_d) = 0 \},\]
where the valuation $v$ is extended to $\overline{K}$.

\section{Bruhat-Tits buildings} 
\subsection{Basic facts} We fix a ground field $K$ which is complete with respect to a non-trivial discrete valuation $v : K^\ast \rightarrow \R$ with perfect residue field $k$. By $\mathcal{O}_K$ we denote the ring of integers in $K$, i.e. $\mathcal{O}_K = \{ x \in K:v(x) \geq 0\}$.  For example, $K$ could be a local field, i.e. locally compact with respect to the topology induced by $v$. In this case, $K$ is either a finite extension of $\mathbb{Q}_p$ or a field of formal Laurent series over a finite field. Alternatively, $K$ could be the field of formal Laurent series $k((T))$ over any perfect ground field $k$.

Let $G$ be a connected reductive group over $K$. In the two groundbreaking papers \cite{brti1} and \cite{brti2}, Bruhat and Tits  define a metric space $\mathfrak{B}(G)$, now called Bruhat-Tits building, endowed with a continuous $G(K)$-action. 

Let $T$ be a maximal split torus in $G$. Then the space $A(T) = X_\ast(T)_\R$ is called the apartment associated to $T$. For every $t \in T(K)$ let $\nu(t)$ be the point in $A(T)$ such that $\chi(\nu(t)) = -v(\chi(t))$ for all $\chi \in X^\ast(T)$. Then $t \in T(K)$ acts on $A(T)$ as the translation with $\nu(t)$.  Let $N$ be the normalizer of $T$ in $G$. Then this action can be continued to an action of $N(K)$ on $A(T)$ by affine-linear maps. 

In \cite{brti2}, Bruhat and Tits construct for every $x \in A(T)$ a subgroup  $P_x$ of $G(K)$, which can be used to define the Bruhat-Tits building $\mathfrak{B}(G)$ as the quotient of $G(K) \times A(T)$ after the following equivalence relation:
\[
(g,x) \sim (h,y),  \begin{array}{l}
\mbox{ if and only if  there exists an element } n \in N(K) \\
 \mbox{ such that } nx= y \mbox{ and } g\inv h n \in P_x.
 \end{array}
 \]
The quotient space $\mathfrak{B}(G)$ is endowed with the product-quotient topology and admits a natural continuous $G(K)$-action via left multiplication in the first factor. For all $x \in A(T)$, the group $P_x$ is the stabilizer of $x$ in $G(K)$. All subsets of $\mathfrak{B}(G)$ of the form $g A$ for $g \in G(K)$ are called apartments. In every apartment there is an infinite arrangement of affine hyperplanes, defining a decomposition of the apartment into faces.  

A parahoric subgroup of $G(K)$ is the connected stabilizer of a face as in \cite{brti2}, Definition 5.2.6. If $G$ is semisimple and simply connected, then by \cite{brti2}, Proposition 4.6.32, the parahoric subgroup associated to the face containing $x$ in its relative interior coincides with $P_x$. Besides, there exists a smooth group scheme $\mathfrak{P}_x$ over $\mathcal{O}_K$ with generic fibre $G$ such that $P_x = \mathfrak{P}_x(\mathcal{O}_K)$. Note that if $x$ lies in the relative interior of a face $F$ with respect to the polysimplicial structure on $A$, then $P_x$ stabilizes the whole face $F$. A parahoric corresponding to a point in the interior of an alcove (i.e. a maximal face) is called an Iwahori subgroup. 

Assume that $G$ is semisimple and simply connected, and consider a point $x$ in $A$ contained in the relative interior of a face $F$. By  $\mbox{Star}(x)$ we denote  the union of all faces in $A$ containing $F$. By \cite{brti2}, Th\'eor\`eme 4.6.33, there exists a bijection between the set of faces in $\mbox{Star}(x)$ and the set of parabolic subgroups of $\mathfrak{P}_x \otimes_{\mathcal{O}_K} k$. This bijection has the following property: If the face $F'$ in $\mbox{Star}(x)$ maps to the parabolic subgroup $Q$ of $\mathfrak{P}_x \otimes_{\mathcal{O}_K} k$, then the preimage of $Q(k)$ via the reduction map $\mathfrak{P}_x (\mathcal{O}_K) \rightarrow (\mathfrak{P}_x \otimes_{\mathcal{O}_K} k) (k)$ is equal to the stabilizer of $F'$ in $G(K)$.

\begin{ex} \rm {\bf The special linear group} 

Let $G = SL_{n,K}$ be the special linear group over $K$, and let $T$ be the maximal split torus of diagonal matrices. 

 We define $a_i \in X^\ast(T)$ by
\[a_i (\diag(s_1, \ldots, s_n)) = s_i,\]
where $\diag(s_1, \ldots, s_n)$ is the diagonal matrix with entries $s_1, \ldots, s_n$. 
Then $X^\ast(T) = \bigoplus_{i=1}^n \Z a_i / \Z (a_1 + \ldots + a_n)$, and the root system of $T$ in $SL_{n,K}$ is equal to
\[\Phi=\{a_{i,j}= a_i / a_j: i \neq  j \mbox{ in } \{1, \ldots,n\}\}.\]
It is of type $A_{n-1}$.  

Besides, let $X_\ast(T) = \Hom_{K} (\Gm, T)$ denote the cocharacter group of $T$. 
Let $\eta_i : \Gm \rightarrow GL_{n,K}$ be the cocharacter of $GL_{n,K}$ mapping $x$ to ${\rm diag}(t_1,\ldots, t_n)$, where $t_i= x$ and $t_j = 1$ for $j \neq i$. Then 
\[X_\ast(T) = \{ m_1 \eta_1 + \ldots + m_n \eta_n : m_i \in \Z \mbox{ with }\sum_i m_i = 0\}.\]

The $\R$-vector space $A = X_\ast(T) \otimes_\Z \R = \{\sum_{i=1}^n x_i \eta_i : x_i \in \R \mbox{ with } \sum_i x_i = 0\}$ is the apartment given by the torus $T$ in the Bruhat-Tits building $\mathfrak{B}(SL_{n,K})$. 
Mapping $\eta_1,\ldots, \eta_n$ to the canonical basis of $\R^n$ provides a homeomorphism
\[A \longrightarrow \T^{n-1}\]
between the apartment $A$ and the tropical torus $\T^{n-1}$. 

For every  $t = \diag(t_1,\ldots, t_n) \in T(K)$  with entries $t_1,\ldots, t_n \in K^\ast$ we define a point in $A$ by $\nu(t) = -v(t_1) \eta_1 + \ldots + -v(t_n) \eta_n$. Then $ t \in T(K)$ acts on $A$ by translation with $\nu(t)$. Besides, let $N$ be the normalizer of $T$ in $SL_{n,K}$. For every element $n \in N(K)$  there is a  permutation $\sigma$ on $\{1, \ldots, n\}$ such that  $n(e_i) = t_{i} e_{\sigma(i)}$ for suitable $t_1, \ldots, t_n \in K^\ast$. The Weyl group $W = N(K) / T(K)$ can therefore be identified with the 
symmetric group on $n$ elements. Hence $W$ acts in a natural way on $A$ by permuting the coordinates of a given point. We can put both actions together to an action of $N(K)$ on $A$ by affine-linear transformations. 

The simplicial structure on $A$ is defined via 
the cells in the infinite hyperplane arrangement consisting of all affine hyperplanes of the form 
\[
H^{(ij)}_{m} \; = \; \biggl\{ \sum_{\ell=1}^n r_\ell \eta_\ell \in A: r_i - r_j = m \biggr\}
\qquad \hbox{
for $ 1 \leq i  < j \leq n$ and $m \in \Z$.}
\]
The isomorphism between $A$ and the tropical torus maps the vertices in $A$ to the subset $\Z^n / \Z (1, \ldots, 1)$ of $\T^{n-1}$.

Let us now describe some of the groups $P_x$. 
For $x = 0$ we have $P_x = SL_n(\mathcal{O}_K)$. If $y = nx \in A$ for some $n \in N(K)$, then $P_y = n SL_n(\mathcal{O}_K) n^{-1}$. For every element $M = (m_{ij})$ of $SL_n(\mathcal{O}_K)$ we  denote by $\overline{M}$ the matrix in $SL_n(k)$ with entries $\overline{m_{ij}}$, where $\overline{m_{ij}}$ is the image of $m_{ij}$ under the residue map $\mathcal{O}_K \rightarrow k$. As explained above, if $y \in A$ is any point in $\mbox{Star}(0)$, then there exists a parabolic subgroup $Q \subset SL_{n,k}$ such that $P_y = \{M \in SL_n(\mathcal{O}_K): \overline{M} \in Q(k)\}$, and every parabolic subgroup gives rise to some $P_y$ in this way.
In other words, for every flag $\mathcal{F}$ of linear subspaces in $k^n$ there exists a point $y \in \mbox{Star}(0)$ such that 
\[P_y = \{ M \in SL_n(\mathcal{O}_K): \, \overline{M} \mbox{ stabilizes } \mathcal{F}\}.\]
For example, the Iwahori group 
\[I = \{ M = (m_{ij}) \in SL_n(\mathcal{O}_K): v(m_{ij}) > 0 \mbox{ for all }i>j\}\]
consisting of  matrices in $SL_n(K)$ with such that all entries  below the diagonal have positive valuation and such that all other entries have non-negative valuation occurs in this way. 

Note moreover that $\mathfrak{B}(PGL_{n,K})$ and $\mathfrak{B}(SL_{n,K})$ are isomorphic, and that $\mathfrak{B}(PGL_{n,K})$ can be identified with the Goldman-Iwahori space of all non-Archimedean norms on $K^n$ modulo scaling, see \cite{goi} and \cite{brti3}. Here a non-Archimedean norm is a map $\gamma: K^n \rightarrow \mathbb{R}_{\geq 0}$ satisfying the following conditions: $\gamma(\lambda v) = |\lambda| \gamma (v)$ and $\gamma (v+w) \leq \sup\{\gamma(v),\gamma(w)\}$ for all $\lambda \in K$ and $v,w \in K^n$,  and $\gamma(v) = 0$ implies $v=0$. 
Via this identification, the apartment $A$ consists of all norms (modulo scaling) of the form
\[\gamma((\lambda_1, \ldots, \lambda_n))= \sup\{|\lambda_1| r_1, \ldots, |\lambda_n| r_n\} \]
for some real vector $(r_1,\ldots, r_n)$. 

There is a dual description in terms of lattices. Namely, the simplicial structure on $\mathfrak{B}(PGL_{n,K})$ and hence on $\mathfrak{B}(SL_{n,K})$ can be described as a flag complex whose vertex set consists of all homothety classes of $\mathcal{O}_K$-lattices in $K^n$. Here two lattice classes are adjacent  if and only if there are representatives $M$ and $N$ of these two classes satisfying $\pi M \subset N \subset M$, where $\pi$ is a prime element in the ring of integers $\mathcal{O}_K$. Let $\{[M_1], \ldots, [M_r]\}$ be a face in the building, i.e. a set of pairwise adjacent lattice classes. We choose the representatives $M_i$ for $i \geq 2$ such that $\pi M_1 \subset M_i \subset M_1$. Then the subspaces $M_i / \pi M_1$ of the $k$-vector space $M_1 / \pi M_1$ form a flag. The stabilizer of this flag over $k$ is precisely the parabolic subgroup $Q$ describing the stabilizer of any point in the interior of the face $\{[M_1], \ldots, [M_r]\}$.

\end{ex}

\begin{ex} \rm {\bf: The symplectic group} 

We consider the semisimple group $G = Sp_{2n,K}$. 
Write $J$ for the $n \times n$ matrix given by 
\[
J = \left(
\begin{array}{lllll} &&&&1\\&0&&1&\\&&\ldots&&\\&1&&0&\\1&&&& \end{array} \right) \] and
set \[ \Psi = \left(\matrix{0&J\cr-J&0}\right)\, .  
\] 
Then $\Psi$ is the coordinate matrix of a standard symplectic form, and
\[Sp_{2n}(K) = \{ M \in SL_{2n}(K): {}^t \!\! M \Psi M = \Psi\}\]
is the group of all matrices preserving the symplectic form given by $\Psi$. 
If $M = (M_{ij})$ is a $n \times n$
matrix, put $
M^\dagger = J\, {}^t\!\!M J $.
Then $M^\dagger$ is obtained from $M$ by  reflection in the anti-diagonal, i.e. 
$M^\dagger_{i,j} = M_{n+1-j,n+1-i}$.
Then 
\[Sp_{2n}(K) = \left\{
\left(\begin{array}{ll} A&B\\C&D \end{array} \right):  A^\dagger D - C^\dagger B = 1,
A^\dagger C = C^\dagger A, B^\dagger D = D^\dagger B \right\}\, . \]
The subgroup $T$ of all diagonal matrices with diagonal entries $(s_1, \ldots, s_n, s_n^{-1}, \ldots, s_1^{-1})$ is a maximal torus in $Sp_{2n,K}$. 
The corresponding root system is $\Phi(T, Sp_{2n,K}) = \{\pm 2 a_i: i = 1,\ldots, n\} \cup \{\pm a_i \pm a_j: 1 \leq i < j \leq n\}$, where $a_i$ maps a diagonal matrix to its  $i$th diagonal entry. It is of type $C_n$. 

The apartment $A(T)$ is the real vector space with basis $\eta_1, \ldots, \eta_n$ satisfying $a_i(\eta_j) = \delta_{ij}$ for $i,j \in \{1,\ldots, n\}$. The affine hyperplanes defining the simplicial structure are $\{2a_i(x) = k\}$ for $i = 1, \ldots,n$ and $k \in \Z$ and  $\{\pm a_i \pm a_j = k\}$ for $i < j$ and $k \in \Z$. 

If $x = 0$, then $P_x = Sp_{2n}(\mathcal{O}_K)$, the group of symplectic matrices over the ring of integers. If $y \in \mbox{Star}(x)$ is a point in the star of $y$, then there exists a parabolic subgroup $Q$ of $Sp_{2n,k}$ such that $P_y = \{M \in Sp_{2n}(\mathcal{O}_K): \overline{M} \in Q(k)\}$, where $\overline{M}$ is the matrix over $k$ induced by $M$. The parabolic subgroups in $SL_{2n,k}$ are the stabilizers of flags of totally isotropic subspaces in $k^{2n}$. For example, the Iwahori group 
\[\{M = (m_{ij}) \in Sp_{2n}(\mathcal{O}_K): v(m_{ij}) > 0 \mbox{ for all }i>j\}\]
arises in this way.   
\end{ex}

\subsection{ A tropical view on stabilizer groups} We will now show that the stabilizer groups $P_x$ defined above with Bruhat-Tits theory have a tropical interpretation via matrices stabilizing points under tropical linear operations. We deal with the $SL_{n,K}$-case first and deduce the general case from it. 

\begin{defi} Let $g$ be an element in $SL_n(K)$, and let $x = ~^t(x_1, \ldots, x_n)$ be a point in $\R^n$. We say that $g$ stabilizes $x$ tropically, if $g_{\trop} \cdot ~^t(x_1, \ldots, x_n) = ~^t(x_1, \ldots, x_n)$ holds. 
\end{defi}
By Definition 2.1, we have
\[g_{\trop} \cdot ~^t(x_1,\ldots, x_n) = ~^t( \max_j \{-v(g_{1j}) + x_j\}, \ldots, \max_j \{-v(g_{nj}) + x_j \}),\]
where $g_{ij}$ are the entries of the matrix $g$. 
Hence $g$ stabilizes $x$ tropically if and only if 
$\max_j \{-v( g_{ij}) + x_j\} = x_i$ for all $i = 1, \ldots, n$. 

\begin{prop} Let $T$ be the torus of diagonal matrices in $SL_{n,K}$ as in Example 2.1, and let $x = \sum_i x_i \eta_i$ be a point in the apartment $A$. Then the stabilizer group $P_x$  with respect to the action of $SL_n(K)$ on $\mathfrak{B}(SL_{n,K})$ is equal to the set of all elements in $SL_n(K)$ stabilizing $~^t(x_1, \ldots, x_n)$ tropically. 
\end{prop}

Recall that the map $(g,x) \mapsto g_{\trop} \cdot x$ does not define an action of $SL_n(K)$ on $\R^n$. Hence the Proposition also shows that the set of matrices in $SL_n(K)$ stabilizing  a point in $\R^n$ tropically is a group, which is not a priori clear. 

{\bf Proof: }First we consider the case that every  $x_i$ lies in the image of the valuation map $v: K^\ast \rightarrow \R$. Hence there exist elements $t_i \in K^\ast $ satisfying $x_i = - v(t_i)$. Then $g$ stabilizes $~^t(x_1,\ldots, x_n)$ tropically if and only if
\[\max_j \{-v( g_{ij} t_j  t_i^{-1})\} = 0 \mbox{ for all  } i.\]
  Now $g_{ij}  t_i^{-1} t_j$ is the entry at position $(i,j)$ of the matrix $ t^{-1} g t$, where $t $ denotes the diagonal matrix with entries $t_1, \ldots, t_n$. Therefore the matrix $t^{-1} g t$ lies in $SL_n(\mathcal{O}_K)$. Besides, every $g$ such that $t^{-1} g t$ lies in $SL_n({\mathcal{O}}_K)$ stabilizes $~^t(x_1,\ldots, x_n)$ tropically. Therefore 
\[\{g \in SL_n(K): g \mbox{ stabilizes } ~^t(x_1, \ldots, x_n) \mbox{ tropically }\} = t SL_n(\mathcal{O}_K) t^{-1}.\]
Since $t SL_n({\mathcal{O}}_K) t^{-1}$ is the stabilizer of $x = \sum_i -v(t_i) \eta_i$, our claim follows.

For a general point $x = \sum_i x_i \eta_i \in A$ there exists a non-Archimedean extension field $L$ of $K$ such that all $x_i$ are contained in the image of the valuation map of $L$ and a continuous $SL_n(K)$-equivariant embedding $\mathfrak{B}(SL_{n,K}) \hookrightarrow \mathfrak{B}(SL_{n,L})$, see e.g. \cite{rtw1}, (1.2.1) and (1.3.4). We have just shown that $\{ g \in SL_n(L): g \mbox{ stabilizes } ~^t(x_1, \ldots, x_n) \mbox{ tropically }\}$ is equal to the stabilizer of $x$ in the building over $L$. Intersecting with $SL_n(K)$ our claim follows.\hfill$\Box$

Now we consider an arbitrary connected semisimple group $G$ over $K$. Let  $\rho: G \rightarrow SL_{n,K}$ be a faithful algebraic representation of $G$, i.e. a homomorphism of $K$-group schemes with trivial kernel.
Let $T$ be a maximal $K$-split torus in $G$, and denote by  $A = A(T)$ the corresponding appartment in $\mathfrak{B}(G)$.  Choose a special vertex $v$ in $A$, i.e. a vertex lying in affine hyperplanes in all possible directions. By \cite{la}, there exists 
 a maximal split torus $T'$ in $SL_{n,K}$ containing $\rho(T)$, and there exists a point $v'$ in the appartment $A' = A(T') = X_\ast(T')_\R$ in $\mathfrak{B}(SL_{n,K})$ given by $T'$
such that the following properties hold:
\begin{enumerate}
\item There is unique affine-linear map $i:A \rightarrow A'$ such that $i(v) = v'$, whose linear part is the map on cocharacter groups 
given by $\rho: T \rightarrow T'$. 

\item The map $i$ satisfies $\rho(P_x) \subset P'_{i(x)}$ for all $x \in A$, where $P_x$ denotes the stabilizer of the point $x$ with respect to the $G(K)$-action on $\mathfrak{B}(G)$, and $P'_{i(x)}$ denotes the stabilizer of the point $i(x)$ with respect to the $SL_n(K)$-action on $\mathfrak{B}(SL_{n,K})$. 

\item Let $Z$ be the centralizer of $T$ in $G$. The map $\rho_\ast: A \rightarrow A' \rightarrow \mathfrak{B}(SL_{n,K})$ defined by composing $i$ with the natural embedding of the appartment $A'$ in the building $\mathfrak{B}(SL_{n,K})$ is $Z(K)$-equivariant, i.e. for all $x \in A$ and $n \in Z(K)$ we have $\rho_\ast (nx) = \rho(n) \rho_\ast(x)$. 

\end{enumerate}

Assume that $T'$ is a maximal torus in $SL_{n,K}$ containing $\rho(T)$ and that $v'$ is a point in $A(T')$ such that properties 1. to 3. are satisfied. 
Then $\rho_\ast: A \rightarrow \mathfrak{B}(SL_{n,K})$ can be continued to a map
$\rho_\ast: \mathfrak{B}(G) \rightarrow \mathfrak{B}(SL_{n,K})$, which is continuous and $G(K)$-equivariant. By \cite{la}, 2.2.9, $\rho_\ast$ is injective and isometrical, if the metric on $\mathfrak{B}(G)$ is normalized correctly. 

The following result states that the subgroup of $G(K)$ stabilizing a point $x$ in the building $\mathfrak{B}(G)$ coincides with the set of elements in $G(K)$ mapping via $\rho$ to matrices stabilizing $i(x)$ tropically. If $G$ is a semisimple, simply connected classical group, and $\rho$ is the  natural embedding of $G$ in $SL_{n,K}$, this result provides a tropical interpretation of the parahoric subgroups of $G$.
 
\begin{thm} let $G$ be a semisimple group over $K$ and let $\rho: G \rightarrow SL_n$ be a faithful algebraic representation. We fix a maximal split torus $T$ in $G$ and a torus $T'$ in $SL_{n,K}$ with the properties 1. to 3. described above. Then for every point $x$ in the apartment $A = A(T)$ of $\mathfrak{B}(G)$ the group $P_x$ of elements in $G(K)$ stabilizing $x$ is equal to 
\[P_x = \{g \in G(K): \rho(g) \mbox{ stabilizes } i(x) \mbox{ tropically}\}.\]
\end{thm}

{\bf Proof: } We claim that $P_x = \{g \in G(K): \rho(g) \in P'_{i(x)}\}$. Since $\rho(P_x) \subset P'_{i(x)}$ by property 3. above, we have one inclusion. Assume that $\rho(g)$ is contained in $P'_{i(x)}$. The point $gx \in \mathfrak{B}(G)$ is mapped to $\rho(g)(i(x)) = i(x) $ via $\rho_\ast: \mathfrak{B}(G) \rightarrow \mathfrak{B}(SL_n)$. Since $i(x) = \rho_\ast(x)$ and $\rho_\ast$ is injective, we find that $g \in P_x$. 
Therefore our claim follows from Proposition 2.3.\hfill$\Box$

\begin{ex} \rm  We consider the natural inclusion $\rho$  of $Sp_{2n,K}$ into $SL_{2n,K}$. Let $T'$ be the torus of diagonal matrices in $SL_{2n,K}$, and let $T$ be the subtorus of diagonal matrices contained in  $Sp_{2n,K}$. We use the notation from Examples 2.1 and 2.2. The corresponding map $i: A(T) \rightarrow A(T')$ maps $\eta_i$ to $\eta_i - \eta_{2n+1-i}$, and hence $i(0) = 0$. For these two vertices the properties 1. to 3. above are satisfied. Hence by Theorem 2.5, we find that for all points $x = \sum_{i=1}^n x_i \eta_i \in A(T)$ the stablilizer group $P_x \subset Sp_{2n}(K)$ is equal to the group of symplectic $2n\times 2n$-matrices stabilizing the real vector $~^t(x_1, \ldots, x_n, -x_n, \ldots, -x_1)$ tropically. Recall that for all points $x$ in the star of $0$, the group $P_x$ can be described explicitely as the group of all elements in $Sp_{2n}(\mathcal{O}_K)$  whose reduction modulo $\pi$ lies in a parabolic subgroup over $k$ associated to the simplicial position of $x$. 

\end{ex}

\section{Compactifications of Bruhat-Tits building}

\subsection{Basic facts and examples}

In \cite{we3} for every connected, semisimple group $G$ over a non-Archimedean local field $K$ and for every faithful, geometrically irreducible algebraic representation $\rho: G \rightarrow GL(W)$  a compactification of $\mathfrak{B}(G)$ is constructed. Its boundary can be identified with the union of Bruhat-Tits buildings associated to certain types of parabolics in $G$. The strategy is the following: We use the combinatorics of the weights given by the representation $\rho$ to define a fan in one apartment. This fan leads to a compactification of the apartment. Then we generalize Bruhat-Tits theory to define groups $P_x$ for all $x$ in the compactified apartment and glue all compactified apartments together as in the definition of $\mathfrak{B}(G)$. 

There is a more general approach. Namely, in \cite{rtw1} we realize the Bruhat-Tits building $\mathfrak{B}(G)$ inside the Berkovich analytic space $G^{an}$ and use the projection to analytical flag varieties of $G$ to obtain a family of compactifications of $\mathfrak{B}(G)$. This fits together with the approach in \cite{we3} by \cite{rtw2}. 

Let us recall some facts from \cite{we3} in more detail and give some examples. In this section we assume that $K$ is a non-Archimedean local field.
Let $G$ be a semisimple group over $K$, and let $\rho: G \rightarrow SL_{n,K}$ be a faithful, geometrically irreducible algebraic representation of $G$. Fix a maximal split torus $T$ in $G$. For every basis $\Delta$  of the root system $\Phi(T,G)$ we denote by $\mu_0(\Delta)$ the corresponding highest weight of $\rho$. 
\begin{defi}
We define the fan $\mathcal{F}_\rho$ in $A(T)$ as the set of all faces of the cones
\[C_\Delta(\rho) = \{ x \in A: \mu_0(\Delta) (x) \geq \mu(x) \mbox{ for all weights }\mu \mbox{ of }\rho\},\] 
where $\Delta$ runs over the bases of $\Phi(T,G)$. 
\end{defi}

Note that for every $\rho$ the cone $C_\Delta(\rho)$ contains the Weyl cone $\mathfrak{C}(\Delta) = \{x \in A: a(x) \geq 0 \mbox{ for all }a \in \Delta\}$. Since the union of all Weyl cones is the total space $A$, we deduce that $\mathcal{F}_\rho$ has support $A$. The Weyl cones for different bases are different. However, we have
$C_\Delta(\rho) = C_{\Delta'}(\rho)$, whenever $\mu_0(\Delta) = \mu_0(\Delta')$, and this may happen for  $\Delta \neq \Delta'$.

 The fan $\mathcal{F}_\rho$ can be used to define a compactification $\overline{A}_\rho$ of $A$, see
\cite{we3}, Section 2. In fact, we put $\overline{A}_\rho = \bigcup_{C \in \mathcal{F}_\rho} A / \langle C \rangle$ and we endow this space with a topology given by tubular neighbourhoods around boundary points. For a more streamlined definition of this topology see \cite{rtw1}, appendix B.

The fan $\mathcal{F}_\rho$ and hence the compactification $\overline{A}_\rho$ only depend on the Weyl chamber face containing the highest weight of $\rho$, see \cite{we3}, Theorem 4.5. Hence we obtain a finite family of compactifications of $A$ in this way.

Using a generalization of Bruhat-Tits theory one can define a subgroup $P_x$ for all $x \in \overline{A}_\rho$ such that for $x \in A$ we retrieve the stabilizer groups in the building, see \cite{we3}, Section 3. Then we define a compactification $\overline{\mathfrak{B}}(G)_\rho$ of $\mathfrak{B}(G)$  as the quotient of 
the topological space $G(K) \times \overline{A}_\rho$ by the equivalence relation
\begin{eqnarray*}
(g,x) \sim (h,y) & \mbox{if and only if  there exists an element } n \in N \\
~ & \mbox{such that } nx= y \mbox{ and } g\inv h n \in P_x.
\end{eqnarray*}

Then $\overline{\mathfrak{B}}(G)_\rho$ is a compact space with a continuous $G(K)$ action, and for every $x \in \overline{A}_\rho$ the group $P_x$ is equal to the group of elements in $G(K)$ stabilizing $x$. 

\begin{ex} \rm Suppose that $G= SL_{n,K}$. We use the notation of Example 2.1. Assume that $\rho = \mbox{id}$. The weights of the identical representation are $\{a_1, \ldots, a_n\}$. For the basis $\Delta =\{ a_{12}, a_{23}, \ldots, a_{n-1 n}\}$ of the root system $\Phi(T,SL_{n,K})$ the highest weight is $\mu_0(\Delta) = a_1$. Hence 
\[C_\Delta(\rho)=\{x \in A: a_1(x) \geq a_i(x) \mbox{ for all }i\}= \{ \sum_{i=1}^n x_i \eta_i \in A: x_1 \geq x_i \mbox{ for all }i\}.\] 
Since the Weyl group (which is isomorphic to the symmetric group on $n$ elements) acts simply transitively on the set of bases for the root system, every maximal cone is of the form
\[\Gamma_k = \{\sum_{i=1}^n x_i \eta_i\in A: x_k \geq x_i \mbox{ for all }i\}\]
for some $k = 1,\ldots, n$. Hence the fan $\mathcal{F}_\rho$ consists of the cones $\Gamma_k$ and of all their faces. 

Let us now describe the compactification $\ov{A} = \ov{A}_{id}$ associated to $\rho = \mbox{id}$. This space is also investigated in \cite{we1}, where it is shown that the correponding compactification of the building $\mathfrak{B}(SL_{n,K})$ contains all homothety classes of free $\mathcal{O}_K$-modules whose rank is strictly smaller than $n$ as vertices on the boundary.  

We write $[n] = \{1, \ldots, n\}$. For every non-empty $I \subset [n]$ 
put $D_I = \cap_{i \in I} \Gamma_i$. 
Let $ \langle D_I \rangle$ be the linear subspace of $A$ generated by $D_I$, and put $A_I = A /  \langle D_{I} \rangle$, and denote the quotient map by $r_I : A \rightarrow A_I$.  Put
\[\ov{A} = \bigcup_{\emptyset \neq I \subset [n]} A_I.\]
The topology on $\ov{A}$ is defined with tubular neighbourhoods around points in $A_I$. To be precise, 
for all open and bounded subsets $U \subset A$ put
\[ C_U^I = \bigcup_{I \subset J \subset [n]} r_J (U + D_I).\]
Then the topology on $\ov{A}$ is the topology with the basis consisting of all $C_U^I$ for non-empty $I \subset [n]$ and of all open bounded subsets $U$ of $A$. 

Let us give a more explicit definition of this space. Put $\R_{- \infty}= \Rinf$ and $\R_{-\infty}^n = (\R \cup \{- \infty\})^n$, and let 
\[(\R_{- \infty}^n)_{I} = \{(x_1, \ldots, x_n) \in \R_{- \infty}^n: x_i =- \infty \mbox{ if and only if } i \notin I \}.\] 
We define an equivalence relation $\sim$  on  $(\R_{- \infty}^n)_{I}$ by $(x_1, \ldots, x_n) \sim (y_1, \ldots, y_n)$ if and only if there exists some $a \in \R$ satisfying $x_i + a = y_i$ for all $i$. Here we put  $- \infty + a = - \infty$ for all $a \in \R$. 
Then the map $A \rightarrow (\R_{- \infty}^n)_{I} / \sim$ which associates to $\sum x_i \eta_i$  the point with coordinates
$x_i$ at $i \in I$ and $- \infty$ at all places $i \notin I$, factors through a bijection $A_I \rightarrow (\R_{- \infty}^n)_{I} /  \sim$.
Hence we get a  bijection 
\[\alpha: \ov{A} \longrightarrow (\R_{-\infty}^n \backslash \{(-\infty, \ldots, - \infty)\}) / \sim .\]

We endow $\R_{- \infty}$ with the topology such that the sets $\{x \in \R: x < b \}$ for $b \in \R$ form a basis of neighbourhoods around $- \infty$, and $(\R_{-\infty})^n$ with the product topology. Then one can check easily that $\alpha$ is a homeomorphism, if we endow the space on the right hand side with the quotient topology. 

This compactification of $A$ could be regarded as a tropical analog of projective space.

\end{ex}

\begin{ex} \rm  Suppose that $\rho$ is a representation of $G = SL_{n,K}$ with highest weight $\mu_0(\Delta) = n a_1 + (n-1) a_2 + \ldots + 2 a_{n-1} + a_n$ for $\Delta = \{a_{12}, \ldots, a_{n-1 n}\}$. Then 
\[C_\Delta(\rho) = \{ x \in A: a_{12} (x) \geq 0, \ldots, a_{n-1 n}(x) \geq 0 \} = \mathfrak{C}(\Delta)\]
is the Weyl cone associated to $\Delta$.  Hence the fan $\mathcal{F}_\rho$ is the Weyl fan consisting of all $\mathfrak{C}(\Delta)$ and their faces.

The corresponding compactification $\overline{\mathfrak{B}}(G)_\rho$ coincides with Landvogt's polyhedral compactification studied in \cite{la}.
\end{ex}

\begin{ex} \rm 
Suppose that $G= Sp_{2n,K}$ and that $\rho: Sp_{2n,K} \rightarrow SL_{2n,K}$ is the canonical embedding. We use the same notation as in Example 2.2. The weights of $\rho$ are the characters $a_1, \ldots, a_n, -a_1, \ldots, -a_n$ of $T$. Hence for the basis $\Delta = \{a_1 - a_2, a_2 - a_3, \ldots, a_{n-1} - a_n, 2a_n\}$ of the root system, the highest weight is $a_1$, and we find
\[C_\Delta(\rho) = \{\sum_{i=1}^n x_i \eta_i \in A(T): x_1 \geq 0 , x_1 \geq \max\{x_2, -x_2\} ,\ldots, x_1 \geq \max\{x_n,-x_n\}\}.\]
The other maximal cones in the fan $\mathcal{F}_\rho$ are the translates of $C_\Delta(\rho)$ under the Weyl group, which is the group of signed permutation of $n$ elements. Hence $\mathcal{F}_\rho$ consists of all faces of the $2n$ cones
\[\Gamma_{k,+} = \{\sum_{i=1}^n x_i \eta_i \in A(T): x_k \geq 0 \mbox{ and } x_k \geq \max\{x_j,-x_j\} \mbox{ for all }j \neq k\}\]
and
\[\Gamma_{k,-}= \{\sum_{i=1}^n x_i \eta_i \in A(T): x_k \leq 0 \mbox{ and } x_k \leq \min\{x_j,-x_j\} \mbox{ for all }j \neq k\}\]
for $k = 1, \ldots, n$. 
\end{ex}

\subsection{Weight polytopes}

As in Section 3.1, let $\rho: G \rightarrow SL_{n,K}$ be a faithful, geometrically irreducible representation of the semisimple group $G$, and 
let $T$ be a maximal split torus in $G$. By  $\Pi_{\rho} \subset X^\ast(T)_\R = X^\ast(T) \otimes_\Z \R$ we denote the weight polytope of $\rho$, i.e. $\Pi_\rho$ is the convex hull of all weights of $\rho$. Recall that we identify $X(T)^\ast_\R$ with the dual space of $A(T)$. 

\begin{prop}
The fan $\mathcal{F}_\rho$ from Definition 3.1 is the normal fan of the weight polytope $\Pi_{\rho}$. 
\end{prop}

{\bf Proof: }By definition, every face $F$ of $\Pi_\rho$ gives rise to a face $\mathcal{N}(F)$ of the normal fan defined by
\[\mathcal{N}(F) = \{x \in A(T): F \subset \mbox{face}_x(\Pi_\rho)\}, \]
where 
\[\mbox{face}_x(\Pi_\rho) = \{p \in \Pi_\rho: p(x) \geq q(x) \mbox{ for all } q \in \Pi_\rho\},\]
see \cite{zi}, Section 7.1.
It suffices to check that the cones of maximal dimension in the normal fan of $\Pi_\rho$ coincide with the cones of maximal dimension in $\mathcal{F}_\rho$. 
Let $\Delta$ be a basis of the root system $\Phi(T,G)$, and let $\mu_0(\Delta)$ be the corresponding highest weight of $\rho$. Then $\{ \mu_0(\Delta)\}$ is a vertex in $\Pi_{\rho}$, giving rise to the following  maximal cone of the normal fan:
 \begin{eqnarray*}
 \lefteqn{\mathcal{N}(\{ \mu_0(\Delta)\})}\\
 & = & \{ x \in A(T): \mu_0(\Delta) (x) \geq \mu(x) \mbox{ for all  weights }\mu \mbox{ of } \rho\}\\
 &  = & C_\Delta(\rho)
 \end{eqnarray*}
 Since every vertex of $\Pi_\rho$ is given by a highest weight, our claim is proven.\hfill$\Box$
 
Now we consider the trace of the representation $\rho$ on $T$, i.e. the morphism 
\[ \mbox{tr}(\rho): T \hookrightarrow G \stackrel{\rho}{\rightarrow} SL_{n,K} \stackrel{\mbox{tr}}{\rightarrow} \mathbb{A}_K^1.\]
Then $\mbox{tr}(\rho)$  is given by a global section of $T$, i.e. by an element $f_\rho \in K[X^\ast(T)]$. 

\begin{prop}
The tropicalization of the hypersurface in $T$ given by the polynomial $f_\rho$ is equal to the codimension one-skeleton of the fan  $\mathcal{F}_\rho$.
\end{prop}

{\bf Proof: } By definition, $f_\rho = \sum_{\mu} \dim(V_\mu) \mu$, where the sum runs over the weights of $\rho$, and where $V_\mu = \{ v \in K^n: \rho(t) v = \mu(t) v \mbox{ for all } t \in T(K)\}$ is the weight space associated to $\mu$. Let us put $c_\mu = \dim(V_\mu)$. Recall from Section 1 that the tropicalisation of the hypersurface given by $f_\rho$ is equal to 
\[\mathcal{T}(f_\rho) = \{x \in A(T): \max_\mu \{-v(c_\mu) + \mu(x)\} \mbox{ is attained at least twice}\}.\]
The complement of $\mathcal{T}(f_\rho)$ in $A(T)$ is the set of all points $x$ such that there exists a weight $\mu$ satisfying $-v(c_\mu) + \mu(x) > -v(c_\lambda) + \lambda(x)$ for all weights $\lambda \neq \mu$ of $\rho$. Now $x$ lies in a cone $C_\Delta(\rho)$ of $\mathcal{F}(\rho)$. Since $\rho$ is geometrically irreducible, the weight spaces of highest weights are one-dimensional, so that $-v(c_{\mu_0(\Delta)}) = 0$. Hence we find that $-v(c_\mu) + \mu(x) \leq -v(c_{\mu_0(\Delta)}) + \mu_0(\Delta)(x)$, which implies $ \mu = \mu_0(\Delta)$. Hence $x$ lies in the interior of $C_\Delta(\rho)$. Therefore the complement of $\mathcal{T}(f_\rho)$ is equal to the union of the interiors of all maximal cones of $\mathcal{F}_\rho$, which implies our claim.\hfill$\Box$

\begin{ex} \rm We look at a geometrically irreducible representation $\rho$ of $SL_{n,K}$ and use the notation from Example 2.1. Assume that the characteristic of $K$ is zero. Then $\rho$ is induced by a representation of $GL_{n,K}$, which can be described with a partition $\lambda_1 \geq  \ldots \geq \lambda_n \geq 0$ of  $\lambda_1 + \ldots + \lambda_n$, see e.g. \cite{gr}. The trace of this representation is given by the Schur polynomial $S_\lambda$ associated to $\lambda = (\lambda_1, \ldots, \lambda_n)$, which is equal to 
\[S_{\lambda}(z_1,\ldots, z_n) = \frac{ \det((z_j^{\lambda_i + n -i})_{i,j = 1, \ldots n})}{\det ((z_j^{n-i})_{i,j = 1, \ldots,n})}.\]

Let $T_0$ be the torus of diagonal matrices in $GL_{n,K}$. Then $X^\ast(T_0) = \bigoplus_i \Z a_i$ and  $X^\ast(T) = X^\ast(T_0) / \Z (\sum_i a_i)$. Besides, $X_\ast(T) \subset X_\ast(T_0) = \bigoplus_i \Z \eta_i$.
Hence the polynomial $S_\lambda$ in $K[X^\ast(T_0)]$ maps to $f_\rho$ under the natural map $K[X^\ast(T_0)] \rightarrow K[X^\ast(T)]$. Therefore the tropicalization of the Schur polynomial $\mathcal{T}({S_\lambda}) \subset X_\ast(T_0)_\R $ satisfies $\mathcal{T}(S_\lambda) \cap X_\ast(T)_\R = \mathcal{T}(f_\rho)$. 
Hence $\mathcal{F}_\rho$ is the fan induced by the tropicalization of the Schur polynomial $S_\lambda$. 

If $\rho$ is the identity representation, the Schur polynomial is equal to $S_{(1,0,\ldots,0)} = x_1 + \ldots + x_n$. In this case, the tropicalization of the Schur polynomial is a tropical hyperplane. 
\end{ex}

\subsection{Stabilizers of boundary points}
We want to derive a description of the  stabilizers of boundary points in the compactifications $\overline{\mathfrak{B}}(G)_\rho$ with tropical linear algebra, thereby generalizing Theorem 2.5. First we look at the case of the identity representation of $SL_{n,K}$. The corresponding fan $\mathcal{F}_{id}$ corresponds to the tropical hyperplane by Example 3.7, and the induced compactification $\overline{A}_{id}$ is homeomorphic to the space $(\R_{-\infty}^n \backslash \{(-\infty, \ldots, -\infty)\}) / \sim $ by Example 3.2. 

\begin{prop} Let $x$ be a point in the compactified apartment $\overline{A}_{id}$ given by coordinates $(x_1, \ldots, x_n) \in \R_{-\infty}^n$. The stabilizer $P_x$ of $x$ with respect to the action of $SL_n(K)$ on the compactified building $\overline{\mathfrak{B}}(SL_{n,K})_{id}$ satisfies
\[P_x = \{ g \in SL_n(K): g \mbox{ stabilizes } ~^t(x_1, \ldots, x_n) \mbox{ tropically}\}.\]
\end{prop}
Here we extend Definition 2.3 in an obvious way to vectors $~^t (x_1,\ldots, x_n)$ in  $\R^n_{-\infty}$ using the rule $-\infty + a = - \infty$ for all $a \in \R_{-\infty}$. 

{\bf Proof: }Assume that $x$ lies in the boundary component $A_I$ for $I \subset [n]$. Then $^t(x_1, \ldots, x_n)$ lies in $(\R_{-\infty}^n)_I$, i.e. $x_i = - \infty$ if and only if $i \notin I$. Hence $g = (g_{ij})$ stabilizes $^t(x_1, \ldots, x_n)$ tropically if and only if
\begin{eqnarray*}
\max_j \{-v( g_{ij}) + x_j \} & =&  x_i  \mbox{ for all  } i \in I \mbox{ and} \\
\max_j \{-v(g_{ij}) + x_j \} &=& -\infty \mbox{ for all } i \notin I.
\end{eqnarray*}
This is equivalent to the fact that $g_{ij} = 0$ for $i \notin I$ and $j \in I$ and that $\max_{j \in I} \{-v( g_{ij}) + x_j \}  =   x_i$  for all $ i \in I$. Let $V_I$ be the subspace of $K^n$ generated by all canonical basis vectors $e_i$ for $i \in I$. Then $g$ stabilizes $~^t(x_1, \ldots, x_n)$ tropically if and only if $g$ restricts to an automorphism $g_I$ of $V_I$ which stabilizes the point with coordinates $x_i$ for $i \in I$ tropically. On the other hand, by the proof of Theorem 5.7 in \cite{we1}, the group $P_x$ consists of all $g \in SL_n(K)$ restricting to an automorphism of $V_I$ which stabilizes $x \in A_I$ with respect to the action of $SL(V_I)$ on the building associated to $SL(V_I)$. Hence our claim follows from Proposition 2.4.\hfill$\Box$

Now let $G$ be an arbitrary semisimple group over $K$, and let $\rho: G \rightarrow SL_{n,K}$ be a geometrically irreducible faithful representation. If $A(T)$ is the apartment in $\mathfrak{B}(G)$ associated to a maximal split torus $T$ in $G$, we have seen in Section 2 that there exists a maximal torus $T'$ in $SL_{n,K}$ with $\rho(T) \subset T'$, and an affine-linear map $i : A(T) \rightarrow A(T')$, which can be extended to a continuous, $G(K)$-equivariant embedding $\rho_\ast: \mathfrak{B}(G) \rightarrow \mathfrak{B}(SL_{n,K})$. It is shown in \cite{rtw2}, Section 5 that this embedding induces an embedding of the compactified building $\overline{\mathfrak{B}}(G)_\rho$ in $\overline{\mathfrak{B}}(SL_{n,K})_{id}$. First of all, by \cite{rtw2}, Lemma 5.1, the fan structures on the apartments fit together, i.e.
the preimage of the fan $\mathcal{F}_{id}$ under the  map of apartments  $i: A(T) \rightarrow A(T')$ is the fan $\mathcal{F}_\rho$. 

This implies that the affine-linear map $i: A(T) \rightarrow A(T')$ has an continuous extension to compactified apartments
\[i: \overline{A(T)}_\rho \rightarrow \overline{A(T')}_{id}.\]
By \cite{rtw2} Lemma 5.2,  for all $x \in \overline{A(T)}_\rho$ we have $\rho(P_x) \subset P'_{i(x)}$, where $P'_{i(x)}$ deontes the stabilizer of $i(x)$ with respect to the action of  $SL_n(K)$ on $\overline{\mathfrak{B}}(SL_{n,K})$. Hence there is a continuous, $G(K)$-equivariant map  
\[\overline{\rho}_\ast: \overline{\mathfrak{B}}(G)_\rho \rightarrow \overline{\mathfrak{B}}(SL_{n,K})_{id}\]
extending $\rho_\ast$.
Moreover, it is shown in \cite{rtw2}, Theorem 5.3 that $\overline{\rho}_\ast$ is a homeomorphism onto the closure of the image of $\rho_\ast$. 

\begin{thm} let $G$ be a semisimple group over $K$ and let $\rho: G \rightarrow SL_{n,K}$ be a geometrically irreducible, faithful algebraic representation. We fix maximal split  tori $T$ in $G$ and  $T'$ in $SL_{n,K}$ as above. Then for every point $x$ in the compactified apartment $\overline{A(T)}_\rho$ the group $P_x$ of elements in $G(K)$ stabilizing $x$ is equal to 
\[P_x = \{g \in G(K): \rho(g) \mbox{ stabilizes } i(x) \mbox{ tropically}\}.\]
\end{thm}
{\bf Proof: }We proceed as in the proof of Theorem 2.5. Let $x$ be a point in $\overline{A(T)}_\rho$. By \cite{rtw2}, Lemma 5.2, we have $\rho(P_x) \subset P'_{i(x)}$. Any $g \in G(K)$ with  $ \rho(g) \in P'_{i(x)} $ satisfies $\overline{\rho}_\ast (gx) = \rho(g) \overline{\rho}_\ast(x) = \overline{\rho}_\ast(x)$. By \cite{rtw2}, Theorem 5.3, the map $\overline{\rho}_\ast$ is injective, which implies $g \in P_x$. Therefore $P_x = \rho^{-1} P'_{i(x)}$, and our claim follows from Proposition 3.8.\hfill$\Box$

\begin{center}
Institut f\"ur Mathematik\\
Goethe-Universit\"at Frankfurt\\
Robert-Mayer-Strasse 8\\
D- 60325 Frankfurt\\
email: werner@math.uni-frankfurt.de
\end{center}

\end{document}